\documentclass[11pt]{amsart}

\usepackage{amssymb}

\evensidemargin\oddsidemargin

\begin{document}
\baselineskip=18pt
\setcounter{page}{1}
    
\newtheorem{Conj}{Conjecture}
\newtheorem{Coj}{Conjecture\!\!}
\newtheorem{Corr}{Corollary}
\newtheorem{Prop}{Proposition}
\newtheorem{Theo}{Theorem\!\!}
\newtheorem{Rqs}{Remarks\!\!}
\newtheorem{Lemm}{Lemma}
\newtheorem{Corru}{Corollary\!\!}
\newtheorem{Corrd}{Corollary 2\!\!}
\newtheorem{Rq}{Remark}

\renewcommand{\theCoj}{}
\renewcommand{\theTheo}{}
\renewcommand{\theCorrd}{}
\renewcommand{\theCorru}{}

\def\a{\alpha}
\def\Aa{{\mathcal A}}
\def\Ss{{\mathcal S}}
\def\b{\beta}
\def\B{{\bf B}} 
\def\C{{\mathcal{C}}} 
\def\CC{{\mathbb{C}}} 
\def\E{{\mathcal{E}}} 
\def\EE{{\mathbb{E}}} 
\def\Da{{\rm D}_\a}
\def\Dea{\Delta_\a}
\def\Dq{{\rm D}_q}
\def\esp{{\mathbb{E}}} 
\def\elaw{\stackrel{d}{=}}
\def\eps{\varepsilon}
\def\F{{\bf F}} 
\def\G{\gamma} 
\def\Fa{F_\a}
\def\Ka{K_\a}
\def\Xa{X_\a} 
\def\ca{c_\a}
\def\fa{f_\a} 
\def\ga{g_\a} 
\def\sa{s_\a} 
\def\HH{{\mathbb{H}}} 
\def\hS{{\hat S}}
\def\hT{{\hat T}}
\def\hX{{\hat X}}
\def\ii{{\rm i}}
\def\K{{\bf K}} 
\def\L{{\mathcal{L}}} 
\def\lb{\lambda}
\def\lacc{\left\{}
\def\lcr{\left[}
\def\lpa{\left(}
\def\lva{\left|}
\def\NN{{\mathbb{N}}} 
\def\pa{p_\a}
\def\Pa{P_\a}
\def\Qa{Q_\a}
\def\Ra{R_\a}
\def\pb{{\mathbb{P}}}
\def\bQa{{\bar Q}_\a}
\def\Qqa{{\widehat \Qa}}
\def\R{{\mathcal{R}}}
\def\rl{{\mathbb{R}}}
\def\racc{\right\}}
\def\rcr{\right]}
\def\rpa{\right)}
\def\rva{\right|}
\def\T{{\bf T}} 
\def\TT{{\rm T}} 
\def\U{{\bf U}} 
\def\Un{{\bf 1}}
\def\ZZ{{\mathbb{Z}}} 

\newcommand{\fin}{\vspace{-0.4cm}
                  \begin{flushright}
                  \mbox{$\Box$}
                  \end{flushright}
                  \noindent}

\title[Total positivity of a Cauchy kernel]{Total positivity of a Cauchy kernel}

\author[Thomas Simon]{Thomas Simon}

\address{Laboratoire Paul Painlev\'e, Universit\'e Lille 1, Cit\'e Scientifique, F-59655 Villeneuve d'Ascq Cedex. 
{\em Email}: {\tt simon@math.univ-lille1.fr}}



\begin{abstract} We study the total positivity of the kernel $1/(x^2 + 2 \cos(\pi\a)xy +y^2).$ The case of infinite order is characterized by an application of Schoenberg's theorem. We then give necessary conditions for the cases of any given finite order with the help of Chebyshev polynomials of the second kind. Sufficient conditions for the finite order cases are also obtained, thanks to Propp's formula for the Izergin-Korepin determinant. As a by-product, we give a partial answer to a question of Karlin on positive stable semi-groups.   
\end{abstract}

\maketitle

\section{Introduction}

Let $I$ be some real interval and $K$ some real kernel defined on $I\times I$. The kernel $K$ is called totally positive of order $n$ (${\rm TP}_n$) if 
$$\det \lcr K(x_i, y_j)\rcr_{1\le i,j\le m}\; \ge \; 0$$ 
for every $m\in \{1,\ldots, n\},\, x_1< \ldots < x_m$ and $y_1< \ldots < y_m.$  If these inequalities hold for all $n$ one says that $K$ is ${\rm TP}_\infty$. The kernel $K$ is called sign-regular of order $n$ (${\rm SR}_n$) if there exists $\{\eps_m\}_{1\le m\le n} \in \{-1,1\}$ such that
$$\eps_m \det \lcr K(x_i, y_j)\rcr_{1\le i,j\le m}\; \ge \; 0$$ 
for every $m\in \{1,\ldots,n\},\, x_1< \ldots < x_m$ and $y_1< \ldots < y_m.$ If these inequalities hold for all $n$ one says that $K$ is ${\rm SR}_\infty$. The above four properties are called strict, with corresponding  notations STP and SSR, when all involved inequalities are strict. We refer to \cite{K} for the classic account on this field and its various connections with analysis, especially Descartes' rule of signs. We also mention the recent monograph \cite{P} for a more linear algebraic point of view and updated references.

A function $f : \rl \to\rl^+$ is called P\'olya frequency of order $n\le \infty$ (${\rm PF}_n$) if the kernel $K(x,y) = f(x-y)$ is ${\rm TP}_n$ on $\rl\times\rl.$ When this kernel is ${\rm STP}_n,$ we will use the notation $f \in {\rm SPF}_n.$ Probability densities belonging to the class ${\rm PF}_\infty$ have been characterized by Schoenberg - see e.g. Theorem 7.3.2 (i) p. 345 in \cite{K} - through the meromorphic extension of their Laplace transform, whose reciprocal is an entire function having the following Hadamard factorization:
\begin{equation}
\label{droite}
\frac{1}{\EE[e^{sX}]} \; =\; e^{-\gamma s^2+\delta s}\prod_{n=0}^\infty (1+a_n s)e^{-a_n s}
\end{equation}
with $X$ the associated random variable, $\gamma \ge 0,\, \delta\in\rl,$ and $\sum a_n^2 < \infty.$ The classical example is the Gaussian density. ${\rm PF}_2$ functions are easily characterized by the log-concavity on their support - see Theorems 4.1.8 and 4.1.9 in \cite{K}. But except for the cases $n=2$ or $n=\infty$ there is no handy criterion for testing the ${\rm PF}_n$ character of a given function resp. the ${\rm TP}_n$ character of a given kernel, and such questions might be difficult. In this paper we consider the kernel
$$\Ka(x,y)\; =\; \frac{1}{x^2 + 2 \cos(\pi\a)xy +y^2}$$
over $(0,+\infty)\times(0, +\infty),$ with $\a\in [0,1).$ Because of its proximity with the standard Cauchy density we may call $\Ka$ a Cauchy kernel, eventhough the denomination Poisson kernel would also be justified. Occurrences of the kernel $\Ka$ in the literature are many and varied. We show the following

\begin{Theo} {\em (a)} One has $$\Ka\,\in{\rm STP}_\infty\,\Leftrightarrow\, \Ka\,\in{\rm SR}_\infty\,\Leftrightarrow\, \alpha\in \{1/2, 1/3, \ldots, 1/n, \ldots, 0\}.$$
\noindent
{\em (b)} For every $n\ge 2,$ one has 
$$\Ka\,\in{\rm TP}_n\,\Leftrightarrow\,\Ka\,\in{\rm SR}_n\,\Rightarrow\,\alpha\in \{1/2, 1/3, \ldots, 1/n\}\;\; {\rm or}\;\; \alpha < 1/n.$$ 
\noindent
{\em (c)} For every $n\ge 2,$ one has $\alpha \le  1/n\wedge 1/(n^2-n-6)_+ \,\Rightarrow\, \Ka\,\in{\rm STP}_n.$
 \end{Theo}

The well-known fact that $K_0$ and $K_{1/2}$ are ${\rm STP}_\infty$ is a direct consequence of explicit classical formul\ae\, for the involved determinants, due respectively to Cauchy and Borchardt - see e.g. (2.7) and (3.9) in \cite{Kr} - and which will be recalled thereafter. The characterization obtained in Part (a) follows without difficulty from Schoenberg's theorem and Gauss' multiplication formula, and will be given in Section 2. The more involved proofs of Part (b) and Part (c) rely both on an analysis of the derivative determinant
$$\Dea^n\; =\; {\det}^{}_n \lcr \frac{\partial^{i+j-2}\Ka}{\partial x^{i-1}\partial y^{j-1}}\rcr$$
where, here and throughout, we set ${\det}^{}_n$ for the determinant of a matrix whose rows and columns are indexed by $(i,j)\in \{1,\ldots, n\}^2.$ To obtain Part (b) we establish a closed expression for $\Dea^n(1,0+)$ in terms of Chebyshev polynomials of the second kind, an expression which is negative whenever $ \alpha\not\in \{1/2, 1/3, \ldots, 1/n\}$ or $\alpha > 1/n.$ The computations are performed in Section 3, in four manners involving respectively Wronskians, Schur functions, rectangular matrices and alternating sign matrices. The observation that these four very different approaches all lead to the same formula was interesting to the author. On the other hand only the last approach, which is based on Propp's formula for the Izergin-Korepin determinant, seems successful to get the more difficult Part (c). This latter result is partial because our upper condition on $\a$ is probably not optimal. We raise the
\begin{Coj} For every $n\ge 2,$ one has  
\begin{equation}
\label{Main}
\alpha <  1/n \;\Rightarrow\; \Ka\,\in{\rm STP}_n.
\end{equation}
\end{Coj}
This would show that the inclusion in Part (b) is actually an equivalence, which is true for $n=2,3$ by Part (c) of the theorem because then $n\ge (n^2-n-6)_+.$  In Section 5 we show (\ref{Main}) for $n=4,5$, and we also give some heuristics reasons supporting the validity of this inclusion for all $n$. The latter seems however to require wild  computations. In Section 5, we state three other conjectures of combinatorial nature whose fulfilment would entail (\ref{Main}). The most natural one is a criterion for the positivity of the generating function
$$f_{n,k}(z)\; =\; \sum_{\tiny \begin{array}{c} A\!\in\!\Aa_n\\ \mu(A)\! =\! k\end{array}} \!\!\!\! z^{\nu(A)}{\bar z}^{\nu(A^Q)}$$
evaluated at a certain complex number $z$, where $\Aa_n$ is the set of alternating sign matrices of size $n$, $ A^Q$ is the anticlockwise quarter-turn rotation of $A,$ $\mu(A)$ is the number of negative entries and $\nu(A)$ the inversion number (see the precise notations below). Whereas we can show this criterion for $k = 0$ or $k = \mu_{\rm max}$ (see below Propositions 3 and 4), unfortunately we cannot do so for all $k$ since there is no sufficiently explicit general formula for $f_{n,k}$. Let us stress that the single evaluation of $f_{n,k}(1)=\sharp\{A\in\Aa_n, \, \mu(A)=k\}$ is a difficult open problem, solved only for certain values of $k$ - see \cite{LG} and the references therein.

The present paper was initially motivated by a question of S. Karlin - see Section 6 below for its precise statement - on the total positivity in space-time of the positive stable semi-group $(t,x)\mapsto p_\a(t,x)$ on $(0,+\infty)\times (0, +\infty),$ which we recall to be defined by the Laplace transform
\begin{equation}
\label{PS}
\int_0^\infty  p_\a(t,x) e^{-\lb x}\, dx\; =\; e^{-t\lambda^\a},\qquad    \lb \ge 0.
\end{equation}
As a simple consequence of Part (b), in Section 6 we show the 
\begin{Corru}  For every $n\ge 2,$ one has 
$$\pa\,\in{\rm SR}_n\,\Rightarrow\,\alpha\in \{1/2, 1/3, \ldots, 1/n\}\;\; {\rm or}\;\; \alpha < 1/n.$$ 
\end{Corru}
The inclusion $\a \le 1/n \,\Rightarrow\, \pa\in{\rm TP}_n$ was proved in \cite{ThS1} for $n=2$ and we believe that it is true in general. This would give a complete answer to Karlin's question and also entail the above conjecture - see Section 6 for an explanation. This will be the matter of further research. 
 
\section{Proof of part {\rm (a)}} 

This easy part of the theorem is a consequence of the following Proposition 1 and Corollary 1. 

\begin{Prop} One has  
$$\Ka\,\in{\rm STP}_\infty\;\Leftrightarrow\; \Ka\,\in{\rm TP}_\infty\;\Leftrightarrow\; \alpha\in \{1/2, 1/3, \ldots, 1/n, \ldots, 0\}.$$ 
\end{Prop}

\proof We begin with the second equivalence. Consider the generalized logistic distribution with density 
$$\ga (x)\; =\; \frac{\sin(\pi\a)}{2\pi\a(\cosh (x) + \cos (\pi\a))}$$
over $\rl.$ Simple transformations - see Theorem 1.2.1. in \cite{K} - entail that the ${\rm TP}_\infty$ character of $\Ka$ amounts to the fact that $\ga\in {\rm PF}_\infty.$ For every $s\in (-1,1),$ compute
\begin{eqnarray*}
\int_\rl e^{sx} \ga(x)\, dx & = & \frac{\sin(\pi\a)}{\pi\a}\int_0^\infty \frac{u^s}{u^2 + 2u\cos (\pi\a) + 1} du\; =\; \frac{\sin(\pi\a s)}{\a\sin(\pi s)} 
\end{eqnarray*}
where the right-hand side follows from the residue theorem and is meant as a limit for $\a = 0.$ If $\alpha\notin \{1/2, 1/3, \ldots, 1/n, \ldots, 0\},$ the function
$$s\; \mapsto\; \frac{\a\sin(\pi s)}{\sin(\pi\a s)}$$
has a pole at $1/\a$ so that $\ga\not\in {\rm PF}_\infty$ by the aforementioned Theorem 7.3.2 (i) in \cite{K}. If $\alpha =1/n$ for some $n\ge 2,$ writing 
\begin{equation}
\label{Gauss}
\frac{\sin(\pi s)}{n\sin(\pi s/n)}\; =\; \frac{\Gamma(1-s)\Gamma(1+s)}{\Gamma(1-s/n)\Gamma(1+s/n)}
\end{equation}
and applying Gauss' multiplication formula and Weierstrass formula for the Gamma function - see e.g. 1.2(11) p.4 and 1.1(3) p.1  in \cite{E} - shows that this function is of the type (\ref{droite}), in other words that $\ga\in {\rm PF}_\infty.$ If $\alpha =0,$ the same conclusion holds true thanks to the Eulerian formula
$$\frac{\sin(\pi s)}{\pi s}\; =\; \prod_{n\ge 1} \lpa 1- \frac{s^2}{n^2}\rpa\cdot$$
This finishes the proof of the second equivalence. To show the first one, it remains to prove that $\Ka\,\in{\rm STP}_\infty$ whenever $\alpha\in \{1/2, 1/3, \ldots, 1/n, \ldots, 0\}.$ Cauchy's double alternant formula (see e.g. (2.7) in \cite{Kr} or Example 4.3 in \cite{P}) 
$${\det}^{}_n \lcr \frac{1}{x_i^2 + y_j^2}\rcr\; =\; \frac{\prod_{1\le i<j\le n} (y_j^2 -y_i^2)(x_j^2 - x_i^2)}{\prod_{1\le i,j\le n} (x_i^2 + y_j^2)}$$
entails immediately that $K^{}_{1/2}\in {\rm STP}_\infty.$ Analogously, Borchardt's formula (see e.g. (3.9) in \cite{Kr}) 
\begin{equation}
\label{Borch}
{\det}^{}_n \lcr \frac{1}{(x_i + y_j)^2}\rcr\; =\; {\det}^{}_n \lcr \frac{1}{x_i + y_j}\rcr\,\times\,{\rm perm}^{}_n \lcr \frac{1}{x_i + y_j}\rcr
\end{equation}
yields $K^{}_{0}\in {\rm STP}_\infty.$ An alternative way to prove these two latter facts consists in writing
$$\frac{1}{x^2+y^2} \; =\; \int_0^\infty e^{-x^2u} e^{-uy^2} \, du\quad \mbox{and}\quad\frac{1}{(x+y)^2} \; =\; \int_0^\infty e^{-xu} e^{-uy} \,udu.$$
Indeed, the composition formula - see Lemma 3.1.1 in \cite{K} - entails that $K_{0}$ and $K_{1/2}$ are ${\rm STP}_\infty$ because the kernel $e^{-ux}$ is ${\rm SRR}_\infty$ on $(0,+\infty)\times(0, +\infty)$ - see \cite{K} p. 18 for the latter fact and p. 12 for an explanation of the notation ${\rm SRR}_\infty.$ 

The remaining cases $\a\not\in\{ 0, 1/2\}$ are slightly more involved. First, it follows from (\ref{Gauss}) and a fractional moment identification that $f_{1/n}$ is the density of the independent sum
\begin{equation}
\label{LTG}
\log (\Gamma_{1/n})\, + \,\cdots\, + \,\log (\Gamma_{(n-1)/n}) \,-\,  
\log (\Gamma_{1/n})\, - \,\cdots\, - \,\log (\Gamma_{(n-1)/n})
\end{equation}
where $\Gamma_t$ is for every $t > 0$ the random variable with 
density 
$$\frac{x^{t-1}e^{-x}}{\Gamma(t)}\Un_{\{x>0\}}.$$
Second, it is straightforward that the density of $\log (\Gamma_t)$ is ${\rm SPF}_\infty.$ This property conveys to $f_{1/n}$ thanks to the above factorization and the composition formula.

\endproof

\begin{Rq} {\em (a) The above argument based on the composition formula yields the ${\rm STP}_\infty$ character of all kernels $(ax+by+c)^{-d}$ for $a,b,d > 0$ and $c\ge 0$ (this is the main result of \cite{DM} -  see Theorem 3.1 therein), in writing
$$\frac{1}{(ax+by+c)^d} \; =\; \int_0^\infty e^{-axu} e^{-byu} u^{d-1}e^{-cu}\, du.$$
(b) It is easy to see that $x\mapsto \Ka(\sqrt{x}, \sqrt{y})$ is a completely monotone function for every $y > 0,$ and that there exists a certain positive finite kernel $L_\a(x,y)$ on $(0,+\infty)\times (0,+\infty)$ such that
$$\Ka(x, y) \; =\; \int_0^\infty e^{-x^2 u}\,L_\a(u,y) du.$$ 
When $\a$ is the reciprocal of an integer, by (\ref{LTG}) it is possible to express $L_\a$ as a convolution of weighted Laplace transformation kernels and show that it is ${\rm RR}_\infty.$ The kernel $L_\a$ is less explicit in the other cases but we feel that its sign-regularity index matches the total positivity index of $\Ka,$ in other terms that 
$$\Ka\,\in{\rm TP}_n\Leftrightarrow L_\a\,\in{\rm RR}_n$$ 
for all $n.$  

}
\end{Rq}

\begin{Prop} The function $\ga$ is positive-definite for every $\a \in [0,1).$
\end{Prop}

\proof The beginning of the proof of Proposition 1 entails by analytic continuation that
$$\int_\rl e^{\ii sx} \ga(x)\, dx \; =\; \frac{\sinh(\pi\a s)}{\a\sinh(\pi s)}$$
for every $s\in \rl.$ Since $\a \in [0,1),$ we can apply the Fourier inversion formula to obtain
$$\ga(x)\; =\; \int_\rl e^{\ii sx} \frac{\sinh(\pi\a s)}{2\pi \a\sinh(\pi s)}\, ds, \qquad x\in\rl.$$
The conclusion follows from Bochner's theorem.

\endproof
\begin{Rq} {\em The above proposition is very well-known, but we gave a proof for the reader's comfort. It is also true that  $\ga^r$ is positive-definite for every $\a \in [0,1)$ and $r > 0$ - see e.g. Exercise 5.6.22 (ii) in \cite{Bh}.}
\end{Rq}
 
\begin{Corr} For every $\a \in [0,1)$ and $n\ge 2,$ one has 
$$\Ka\,\in{\rm TP}_n\;\Leftrightarrow\; \Ka\,\in{\rm SR}_n.$$ 
\end{Corr}

\proof 

Proposition 2 entails that for every $\a\in [0,1)$ the function
$$x\; \mapsto \; \frac{1+\cos(\pi\a)}{\cosh (x) + \cos(\pi\a)}$$
is the characteristic function of the random variable $\Xa$ with density
$$t\; \mapsto\;\frac{(1+\cos(\pi\a))\sinh(\pi\a t)}{\sin (\pi\a) \sinh(\pi t)}\cdot$$
For every $n\ge 2, s_1, \ldots, s_n \in \rl$ and $0<x_1<\ldots < x_n,$ this yields 
$$\sum_{1\le i,j\le n} s_i s_j \Ka(x_i, x_j)\; =\; \frac{1}{2(1+\cos(\pi\a))}\;\esp\lcr\lva\sum_{k=1}^n  t_k e^{\ii y_k \Xa}\rva^2\rcr\; > \; 0,$$
with the notation $t_i = s_i/x_i$ and $y_i = \log (x_i).$ Hence, the quadratic form $[\Ka(x_i,x_j)]_{1\le i,j \le n}$ is positive definite for every $n\ge 2$ and $0<x_1<\ldots < x_n.$ In particular, one has 
$${\det}^{} \lcr \Ka(x_i,x_j)\rcr_{1\le i,j \le n}\; >\; \ 0,$$ 
for every $n\ge 2$ and $0<x_1<\ldots < x_n.$ This shows the required equivalence.

\endproof

\section{Proof of part {\rm (b)}} 

Our argument for this part is the same as the one we have used in \cite{ThS2}, in the framework of confluent hypergeometric functions. Let $X = \{x_1< \ldots< x_n\}$ and $Y= \{y_1< \ldots< y_n\}$ be two sets of positive variables and 
$$V^{}_X \; = \prod_{1\le i<j\le n}\!\!\!\!(x_j - x_i) \qquad \mbox{and}\qquad V^{}_Y \; =\prod_{1\le i<j\le n}\!\!\!\! (y_j -y_i)$$
be the usual Vandermonde determinants. Consider $\Da^n(X,Y) = {\det}^{}_n \lcr \Ka(x_i,y_j)\rcr$
and the aforementioned derivative determinant
$$\Dea^n(x,y)\; =\; {\det}^{}_n \lcr \frac{\partial^{i+j-2}\Ka}{\partial x^{i-1}\partial y^{j-1}}(x,y)\rcr, \quad x, y >0.$$
Using repeatedly the formula
$$\sum_{k=0}^p (-1)^{p-k} C^k_p f(z+k\eps)\; \sim \; f^{(p)}(z)\eps^p, \quad \eps\to 0+$$
which is valid for any smooth real function $f$, elementary operations on rows and columns show that
\begin{equation}
\label{derive}
\Dea^n(x,y)\; =\; {\rm sf}(n-1)^2 \lim_{\eps\to 0+} \frac{\Da^n(X_\eps, Y_\eps)}{V_{X_\eps}V_{Y_\eps}}
\end{equation}
where we have set $x_i^\eps = x +(i-1)\eps, y_i^\eps = y + (i-1)\eps$ for $i = 1\ldots n,$ $X_\eps = (x_1^\eps, \ldots, x_n^\eps), \, Y_\eps = (y_1^\eps, \ldots, y_n^\eps),$ and 
$${\rm sf}(k)\; =\; \prod_{i=0}^k i!$$
for the superfactorial number. By Proposition 2, this entails that $\Dea^n(x,x) \ge 0$ for any $x > 0$ and below it will be established that the inequality is actually everywhere strict - see Remark 5. 

On the other hand, if $\Dea^k(1,0+) < 0$ for some $k\in \{2, \ldots, n\},$ then (\ref{derive}) entails that $\Ka$ is not ${\rm TP}_n$ and hence not ${\rm SR}_n$ by Corollary 1. We will prove that this is the case as soon as $\alpha \ge 1/n$ and $\a\notin\{1/2, 1/3, \ldots, 1/n\}.$ More precisely, setting
$$U^\a_k \; = \; \frac{\sin k\pi\a}{\sin\pi\a}\qquad\mbox{and}\qquad V^\a_n \; =\; \prod_{k=1}^n U^\a_k$$
for every $k,n\ge 1,$  we will show that
\begin{equation}
\label{Maindirect}
\Dea^n(1,0+) \; =\; {\rm sf}(n-1)^2 V^\a_n.
\end{equation}
We give four different proofs of (\ref{Maindirect}), each corresponding to a specific approach to evaluate the derivative determinant $\Dea^n (x,y)$ in closed form. The first two proofs involve classical tools, respectively Wronskians and Schur functions. The last two proofs rely on more elaborate and recent results on the Izergin-Korepin determinant, using combinatorial formul\ae\, due respectively to Lascoux and Propp. The latter formula, involving  alternating sign matrices, will be also used to obtain Part (c). 

\subsection{Wronskians} We fix $\a \in [0,1),$ set $\ca = \cos(\pi\a)$ and $\sa = \sin(\pi\a).$  We use the notations 
$$\Qa(x,y) = y + x \cos(\pi\a),\qquad \Qqa(x,y) = \Qa(y,x),\qquad \Ra(x,y) = \Qa(x,y)\Qqa (x,y),$$ 
and $\Pa(x,y) = x^2 + 2xy \cos(\pi\a) + y^2.$ In the following, we will  often skip the dependence in $(x,y)$ for concision. Observe that
\begin{equation}
\label{basic}
\Ra - \ca\Pa = \sa^2 xy \qquad \mbox{and}\qquad {\widehat Q}_\a^2 - \Pa = -\sa^2y^2.
\end{equation}
Introduce the kernels
$$T_{\a,2p}\; =\; (2p)!\sum_{k=0}^p (-1)^{p-k} C^{p-k}_{p+k} (4\Qa^2)^k\Ka^{p+k+1}$$
and
$$T_{\a,2p+1}\; =\; 2(2p+1)!\,\Qa\,\sum_{k=0}^p (-1)^{p+1-k} C^{p-k}_{p+1+k} (4\Qa^2)^k\Ka^{p+k+2}$$ 
for any integer $p.$ 

\begin{Lemm} For every $n\ge 0,$ one has
$$\frac{\partial^n \Ka}{\partial y^n}\; =\; T_{\a, n}.$$
\end{Lemm} 

\proof The formula is obviously true for $n =0,1.$ By induction, it suffices to show that
$$\frac{\partial T_{\a, 2p}}{\partial y}\; =\; T_{\a, 2p+1}\qquad\mbox{and}\qquad\frac{\partial T_{\a, 2p+1}}{\partial y}\; =\; T_{\a, 2p+2}$$
for any integer $p.$ The latter are elementary computations, whose details are left to the reader.

\endproof

Lemma 1 allows to express $\Dea^n$ as the Wronskian of $(T_{\a, 0}, \ldots, T_{\a, n-1})$ built on $x-$derivatives, in other words one has
\begin{equation}
\label{Wronski}
\Dea^n\; =\; {\det}_n^{}\lcr \frac{\partial^{i-1}T_{\a, j-1}}{\partial x^{i-1}}\rcr.
\end{equation}
Observe that
$$T_{\a,2p}(1,0+)\; =\; (2p)!\sum_{k=0}^p (-1)^k C^{k}_{2p-k} (2\cos(\pi\a))^{2p-2k}\; =\; (2p)!\,\frac{\sin((2p+1)\pi\a)}{\sin (\pi\a)}$$
where the second equality comes after identifying two standard definitions of the Chebyshev polynomial of the second kind $U_{2p}.$ Making the same identification for $T_{\a,2p+1}(1,0+)$ and putting the two together entail 
\begin{equation}
\label{Cheb}
T_{\a,r}(1,0+)\; =\; (-1)^r r!\,\frac{\sin((r+1)\pi\a)}{\sin (\pi\a)}
\end{equation}
for every $r\ge 0.$ We next compute the successive $x-$derivatives of $T_{\a,r}$ at $(1,0+).$ The outline of the proof is analogous to Lemma 1, but it is more involved and so we give the details.

\begin{Lemm} For every $j\ge 1$ and $r\ge 0,$ one has
\begin{equation}
\label{abteil}
\frac{\partial^j T_{\a,r}}{\partial x^j}(1,0+)\; =\; (-1)^j\frac{(r+j+1)!}{(r+1)!}\, T_{\a,r}(1,0+).
\end{equation}
\end{Lemm}

\proof We first consider the case $r = 2p.$ Setting
$$a_{j,r}\; =\; \frac{(r+j+1)!}{(r+1)!}\;,\; E^\a_{j,p} \; =\; \frac{(-1)^j}{(2p)!}\frac{\partial^j T_{\a,2p}}{\partial x^j} \; ,\;F^\a_{k,p}\; =\; (-1)^{p-k} C^{p-k}_{p+k} 4^k \Qa^{2k-1},\; G^\a_{k,p}\; =\; \Qa F^\a_{k,p},$$
we will show by induction that for every $j\ge 1,$ one has
\begin{equation}
\label{indic}
E^\a_{j,p} \; =\;a_{j,r}{\widehat Q}_\a^j\sum_{k=0}^p G^\a_{k,p}\Ka^{p+k+j+1}\; +\; y H^\a_{j,p}
\end{equation}
for some smooth function $H^\a_{j,p}$ on $(0,\infty)\times[0,\infty).$ Specifying to $(x,y) = (1,0+)$ clearly entails (\ref{abteil}). We begin with the case $j=1,$ in computing
\begin{eqnarray*}
E^\a_{1,p} & = & 2\lpa\sum_{k=0}^p (p+k+1) F^\a_{k,p}\Ra\Ka^{p+k+2} \; -\; \ca\sum_{k=1}^p k F^\a_{k,p}\Pa \Ka^{p+k+2}\rpa\\
& = & (2p+2){\widehat Q}_\a\sum_{k=0}^p G^\a_{k,p}\Ka^{p+k+2}\; +\; 2xy\sa^2\sum_{k=1}^p  F^\a_{k,p}\Ka^{p+k+2}\\
& = & (2p+2){\widehat Q}_\a\sum_{k=0}^p G^\a_{k,p}\Ka^{p+k+2}\; +\; y H^\a_{1,p}
\end{eqnarray*}
for some smooth function $H^\a_{1,p}$ on $(0,\infty)\times[0,\infty),$  where in the second equality we used the first identity in (\ref{basic}). Suppose now that (\ref{indic}) holds for some $j\ge 1.$ Differentiating, we obtain
\begin{eqnarray*}
E^\a_{j+1,p} & = & a_{j,r}{\widehat Q}_\a^{j-1}\sum_{k=0}^p (2(p+k+j+1)G^\a_{k,p}{\widehat Q}_\a^2 - 2k \ca\Pa F^\a_{k,p}{\widehat Q}_\a - jG^\a_{k,p})\Ka^{p+k+j+2}\; -\; y \frac{\partial H^\a_{j,p}}{\partial x}\\
& = & (2p+j+2)a_{j,r}{\widehat Q}_\a^{j+1}\sum_{k=0}^p G^\a_{k,p}\Ka^{p+k+j+2} \; +\; 2(\Ra- \ca\Pa){\widehat Q}_\a^j\sum_{k=0}^p k F^\a_{k,p} \Ka^{p+k+j+2}\\
& & \;\; \;\; \;\; \;\; \;\; \;\; \;\; \;\; +\;\; ({\widehat Q}_\a^2- \Pa){\widehat Q}_\a^{j-1}\sum_{k=0}^p G^\a_{k,p}\Ka^{p+k+j+2}\; -\; y \frac{\partial H^\a_{j,p}}{\partial x}\\
& = & a_{j+1,r}{\widehat Q}_\a^{j+1}\sum_{k=0}^p G^\a_{k,p}\Ka^{p+k+j+2}\; +\; y H^\a_{j+1,p}
\end{eqnarray*}
for some smooth function $H^\a_{j+1,p}$ on $(0,\infty)\times[0,\infty),$  where in the second equality we used the two identities in (\ref{basic}). This completes the proof of (\ref{indic}) for $r=2p.$ An analogous induction shows that
$$\frac{(-1)^j}{(2p+1)!}\frac{\partial^j T_{\a,2p+1}}{\partial x^j}\; =\;2\,a_{j,2p+1}\,\Qa{\widehat Q}_\a^j\sum_{k=0}^p (-1)^{p+1-k} C^{p-k}_{p+1+k} (4\Qa^2)^k\Ka^{p+k+j+2}\; +\; y I^\a_{j,p}$$
for every $j\ge 1, p\ge 0,$ where $I^\a_{j,p}$ is a smooth function  on $(0,\infty)\times[0,\infty).$ This yields (\ref{abteil}) for $r=2p+1,$ and finishes the proof. 

\endproof

We finally deduce from (\ref{Wronski}), (\ref{Cheb}), Lemma 2, and Leibniz's formula for determinants, that
\begin{eqnarray*}
\Dea^n(1,0+) & = &  {\rm sf}(n-1)\; V^\a_n\; {\det}^{}_n\lcr\frac{(i+j-1)!}{j!} \rcr\\ & = &  {\rm sf}(n-1)^2\; V^\a_n\; {\det}^{}_n\lcr C^j_{i+j-1}\rcr\;\;= \;\;  {\rm sf}(n-1)^2\; V^\a_n,
\end{eqnarray*}
where the last equality follows from the standard evaluation of a binomial determinant which is left to the reader. This completes the proof of (\ref{Maindirect}). 

\begin{Rq} {\em It does not seem that (\ref{Wronski}) and standard  Wronskian transformations can provide any information which is enough explicit to characterize the everywhere positivity of $(x,y)\mapsto \Dea^n (x,y).$ The latter is crucial for (\ref{Main}) - see Section 4 below.}
\end{Rq}

\subsection{Schur functions} This second approach relies on the following well-known expression for the generating function of Chebyshev polynomials of the second kind:
$$ \sum_{k\ge 0} U^\a_{k+1} z^k\; =\;\Ka(1,-z), \qquad \vert z\vert < 1.$$
This entails 
$$\Da^n(X,Y)\; =\; \frac{1}{{\displaystyle \prod_{j=1}^n x_j^2}}\;{\det}^{}_n\lcr \sum_{k\ge 0} U^\a_{k+1} (y_i z_j)^k  \rcr $$
with the above notation for $X$ and $Y,$ having written $z_j = -x_j^{-1}$ and assuming, here and throughout, $\vert y_i z_j\vert < 1$ for all $i, j = 1,\ldots , n.$ The determinant on the right-hand side, denoted by $\Fa^n(X,Y),$  can be further evaluated by multilinear expansion:
$$\Fa^n(X,Y)\; =\; \sum_{k_1 \ldots k_n = 0}^\infty{\det}^{}_n\lcr z_j^{k_i} \rcr \lpa\prod_{l=1}^n U^\a_{k_l+1} y_l^{k_l}\rpa.$$
Notice that in the multiple sum, all indices $k_i$ can be chosen distinct since otherwise the summand is zero. Setting $\Ss_n$ for the symmetric group of size $n$ and 
$$U_{{\rm K}}^\a \; =\; \prod_{l=1}^n U^\a_{k_l+1}$$
for any $n-$tuple ${\rm K} = (k_1, \ldots, k_n),$ one obtains
\begin{eqnarray*}
\Fa^n(X,Y) & = & \sum_{k_1> \cdots> k_n \ge 0} U_{{\rm K}}^\a\lpa \sum_{\sigma\in \Ss_n} {\det}^{}_n\lcr z_j^{k_{\sigma(i)}} \rcr \lpa\prod_{l=1}^n y_l^{k_{\sigma(l)}}\rpa\rpa\\
& = & \sum_{k_1> \cdots> k_n \ge 0} U_{{\rm K}}^\a\;  {\det}^{}_n\lcr z_j^{k_i}\rcr  {\det}^{}_n\lcr  y_j^{k_i}\rcr,
\end{eqnarray*}
where in the second equality we used twice Leibniz's formula. The above can now be rewritten in terms of Schur functions and Vandermonde determinants. Setting $k_i = n-i +\lambda_i,$ one has
\begin{eqnarray*}
\Fa^n(X,Y) & = & \sum_{\lambda_1\ge\cdots\ge \lambda_n \ge 0} U_{{\rm K}}^\a\;  {\det}^{}_n\lcr z_j^{n-i+\lambda_i}\rcr  {\det}^{}_n\lcr  y_j^{n-i +\lambda_i}\rcr\;
= \; V^{}_Y V^{}_Z\; \sum_{\lambda} U_{{\rm K}}^\a\; s_\lambda (Y) s_\lambda (Z)
\end{eqnarray*}
where in the second equality the sum is meant on all partitions of size $n$ and we use the standard notation for Schur functions, displayed e.g. in Chapter 4 p.124 of \cite{B}. Putting everything together and differentiating with (\ref{derive}) yield after some simplifications
$$\Dea^n(x,y)\; =\; \frac{{\rm sf} (n-1)^2}{x^{n(n+1)/2}} \sum_{\lambda}\lpa\prod_{i=1}^n U^\a_{n-i+\lambda_i+1}\rpa s_\lambda (1,\ldots,1)^2 (-yx^{-1})^{\vert\lambda\vert}$$ 
for all $0<y<x,$ with the notation $\vert \lambda\vert = \lambda_1 +\cdots +\lambda_n.$
By the change of variable $\mu_i = \lambda_{n+1-i},$ the classical formula
$$ s_\lambda (1,\ldots,1) \; =\; \frac{\prod_{1\le i < j\le n} (\lambda_i -\lambda_j + j - i)}{{\rm sf} (n-1)}$$
which can be recovered from the Jacobi-Trudi identity - see e.g. Proposition 4.2 in \cite{B} - and a standard binomial determinant evaluation, entails finally
\begin{equation}
\label{Schu}
\Dea^n(x,y)\; =\; \frac{1}{x^{n(n+1)/2}} \sum_{k=0}^\infty (-yx^{-1})^k\!\!\!\!\!\!\!\!\!\!\!\!\sum_{{\tiny \begin{array}{c} \mu_1\le\cdots \le\mu_n\\
\mu_1 +\cdots + \mu_n = k\end{array}}}\!\!\! \lpa\prod_{i=1}^n\; U^\a_{i+\mu_i}\rpa \prod_{1\le i < j\le n} (\mu_j -\mu_i + j - i)^2.
\end{equation}
We can now compute explicitly $\Dea^n(1,0+)$ because the summation is made on the single partition $(0,\ldots, 0):$ we get 
$$\Dea^n(1,0+)\; =\; \prod_{i=1}^n\; U^\a_{i}\;\times\!\!\!\! \prod_{1\le i < j\le n}\!\! (j - i)^2\; =\; {\rm sf (n-1)^2}\; V^\a_n,$$ 
as required by (\ref{Maindirect}).

\begin{Rq} {\em Because of the alternate signs, the above expression (\ref{Schu}) does not seem very helpful either to study the everywhere positivity of $\Dea^n(x,y).$}
\end{Rq}

\subsection{Rectangular matrices} This third approach hinges upon a certain closed expression for the determinant 
$$\Dq(X,Y)\; =\; {\det}^{}_n \lcr \frac{1}{(x_i + y_j)(qx_i + y_j)}\rcr.$$
The latter is called the Izergin-Korepin determinant in the literature, and appears in the context of the six-vertex model - see Chapter 7 in \cite{B}. We use its evaluation in terms of rectangular matrices separating the variables, which is due to Lascoux - see Theorem $q$ in \cite{L}. It is given by
$$\Dq(X,Y)\; =\; \frac{V^{}_X V^{}_Y}{{\rm P}_q(X,Y)}\,{\det}^{}_n \lcr H_X^{}\times E^q_Y\rcr$$
where ${\rm P}_q(X,Y) \; =\;\prod_{1\le i,j\le n} (x_i + y_j)(q x_i + y_j),$
and
$$H^{}_X \; =\; \lcr h_{k-i} (X)\rcr_{1\le i\le n, 1\le k\le 2n-1},\qquad E^q_Y \; =\; \lcr \frac{q^{k-j+1} -q^{j-1}}{q-1}\, e_{n-k+j-1} (Y)\rcr_{1\le k\le 2n-1, 1\le j\le n}$$
are two rectangular matrices involving the complete resp. elementary symmetric functions, which we recall to be defined through the generating functions
\begin{equation}
\label{gene}
\sum_{k\ge 0}h_k(X) t^k\; =\; \prod_{r=1}^n (1-x_r t)^{-1}\qquad\mbox{and}\qquad \sum_{k\ge 0}e_k(Y) t^k\; =\; \prod_{r=1}^n (1+y_r t).
\end{equation}
Setting $q = e^{2\ii \pi \a}, X^q_x = (xq^{-1/2},\ldots,xq^{-1/2}),$ and $Y_y = (y,\ldots, y),$ Lascoux's formula and (\ref{derive}) yield
$$\Dea^n(x,y) \; =\; {\rm sf}(n-1)^2\, q^{-n(n-1)/4}\Ka(x,y)^{n^2}\, {\det}^{}_n \lcr H_{X^q_x}^{}\times E^q_{Y_y}\rcr.$$ 
By (\ref{gene}), we have $h_r(X^q_x) = C^r_{n+r-1} x^r q^{-r/2},$ whence
$$H_{X^q_x}^{}\; =\; \lcr C^{k-i}_{n+k-1-i}q^{(i-k)/2} x^{k-i}\rcr_{1\le i\le n, 1\le k\le 2n-1}.$$
On the other hand, (\ref{gene}) entails $e_r(Y_y) = C^r_n y^r,$ so that after some simplifications
$$E_{Y_y}^{}\; =\; \lcr C^{n-k+j-1}_{n}q^{(k-1)/2} U_{k+2-2j}^\a y^{n-k+j-1}\rcr_{1\le k\le 2n-1, 1\le j\le n},$$
with our above notation for $U^\a_r.$ 
The $i$-th row of the product $H_{X^q_x}^{}\times E_{Y_y}^{}$ having a factor $q^{(i-1)/2}$, we finally obtain
\begin{equation}
\label{separate}
\Dea^n(x,y) \; =\; {\rm sf}(n-1)^2\, \Ka(x,y)^{n^2}\,{\det}^{}_n \lcr A_n(x) \times B_n^\a(y)\rcr
\end{equation}
with the notations
$$A_n(x)\; =\; \lcr C^{k-i}_{n+k-1-i} x^{k-i}\rcr_{1\le i\le n, 1\le k\le 2n-1}$$
and
$$B_n^\a (y)\; =\; \lcr C^{n-k+j-1}_{n} U_{k+2-2j}^\a y^{n-k+j-1}\rcr_{1\le k\le 2n-1, 1\le j\le n}.$$
Setting $L_n^\a(x,y) = {\det}^{}_n \lcr A_n(x) \times B_n^\a(y)\rcr,$ the Cauchy-Binet formula entails
\begin{equation}
\label{CaBi}
L_n^\a(x,y)\; =\sum_{1\le \sigma_1 < \ldots < \sigma_n \le 2n-1}\!\!\!\!\! {\rm A}_\sigma(x) {\rm B}_\sigma^\a(y)
\end{equation}
where $ {\rm A}_\sigma(x)$ is the $n\times n$ minor obtained from the columns $\sigma_1, \ldots, \sigma_n$ in $A_n(x)$ and ${\rm B}_\sigma^\a(y)$ is the $n\times n$ minor obtained from the rows $\sigma_1, \ldots, \sigma_n$ in $B_n^\a(y).$ By Leibniz's formula, one has
$$ {\rm A}_\sigma(x)\, =\,   {\rm A}_\sigma(1) x^{n_\sigma}\quad\mbox{and}\quad {\rm B}_\sigma^\a(y)\, =\, {\rm B}_\sigma^\a(1) y^{n(n-1) -n_\sigma}$$
with the notation 
$$n_\sigma\; =\; \sum_{i=1}^n (\sigma_i -i).$$
This shows that $L_n^\a(x,y)$ is a homogeneous polynomial of degree $n(n-1)$ with coefficient 
\begin{equation}
\label{coeff}
\sum_{n_\sigma = k} {\rm A}_\sigma(1) {\rm B}_\sigma^\a(1)
\end{equation}
for the term $x^k y^{n(n-1) -k}.$ Besides, by (\ref{separate}) and symmetry, we have $L_n^\a(x,y)=L_n^\a(y,x)$ so that these coefficients are palindromic. We compute
$$L_n^\a(1,0+)\; =\; L_n^\a(0+,1)\; =\;  {\rm A}_{\tilde {\rm Id}}(1) {\rm B}_{\tilde {\rm Id}}^\a(1)$$ 
 with the notation ${\tilde  {\rm Id}} = (n, \ldots, 2n-1).$  One finds immediately ${\rm B}_{\tilde {\rm Id}}^\a(1) = U_1^\a\times\cdots\times U_n^\a$ and some elementary linear transformations yield ${\rm A}_{\tilde {\rm Id}}^\a(1) = 1.$ Since $\Ka(0+,1) = 1,$ we finally deduce (\ref{Maindirect}).
 
\subsection{Alternating sign matrices} In this last approach we use Izergin-Korepin's original formula for the determinant $\Dq(X,Y)$ and its expression in terms of alternating sign matrices, due to Propp, which is given in Exercise 7.2.13 p. 244 in \cite{B}. It reads
\begin{equation}
\label{ASM}
\Dq(X,Y) \; = \; \frac{V^{}_X V^{}_Y}{{\rm P}_q(X,Y)}\sum_{A\in\Aa_n} (-1)^{\mu(A)}(1-q)^{2\mu(A)}q^{C^2_n-I(A)}\prod_{i=1}^n x_i^{\mu_i(A)}y_i^{\mu^i(A)}\!\!\!\!\!\prod_{\tiny \begin{array}{c} 1\!\le\! i,j\! \le\! n \\
a_{ij} = 0
\end{array}
}\!\!\!\!\!\!\!\! (\a_{ij}x_i +y_j)
\end{equation}
where $\Aa_n$ stands for the set of $n\times n$ alternating sign matrices (ASM) viz. those matrices made out of $0$s, $1$s and $-1$s for which the sum of the entries in each row and each column is 1, and the non-zero entries in each row and column alternate in sign. We refer to \cite{B} for a comprehensive account on this topic, and also to Section 3 in the recent paper \cite{Be} for updated results. In (\ref{ASM}), the following notations are used: $\mu_i(A)$ resp. $\mu^i(A)$ is the number of $-1$s in the $i$-th row resp. $i$-th column of $A$, $\mu(A)$ the total number of $-1$s in $A$, $I(A)$ the generalized inversion number of $A$ viz.
$$I(A)\; = \; \sum_{i<k, l<j}a_{ij}a_{kl},$$
and 
$$\a_{ij}\; =\; \lacc\begin{array}{ll} q & \mbox{if}\; \sum_{k\le i}a_{kj} \; =\; \sum_{l\le j}a_{il}, \\
1 & \mbox{otherwise.}
\end{array}\right.$$
Notice that ASM matrices without $-1$ are permutation matrices. In particular, one recovers Borchardt's formula (\ref{Borch}) from (\ref{ASM}) in setting $q = 1$ - see Exercise 7.2.14 in \cite{B}. In the following we will set $J(A) = \sharp\{ (i,j),\; \sum_{k\le i}a_{kj} \; =\; \sum_{l\le j}a_{il}\},$ which is the number of southwest or northeast molecules in the terminology of the six-vertex model, and satisfies the formula
\begin{equation}
\label{SWNE}
J(A)\; =\; 2 I(A)\; -\; 2 \mu(A)
\end{equation} 
(see Exercise 7.1.8 in \cite{B}). Setting $q = e^{2\ii\pi\a},$ we see that (\ref{derive}) entails after some simplifications hinging upon (\ref{SWNE}) the following expression for the derivative determinant $\Dea^n(x,y)$:
\begin{equation}
\label{IKP}
{\rm sf} (n-1)^2\Ka(x,y)^{n^2}\sum_{A\in\Aa_n} (\Pa(x,y))^{\mu(A)} (\Qa(x,y))^{J(A)}(\bQa(x,y))^{n(n-1) - 2I(A)}
\end{equation}
with the notations
$$\Pa(x,y)\; =\; 4\sin^2(\pi\a) xy\qquad\mbox{and}\qquad \Qa(x,y)\; =\; e^{\ii\pi\a/2} x + e^{-\ii\pi\a/2} y.$$
This yields
$$\Dea^n(1,0+)\; =\; {\rm sf} (n-1)^2 e^{-\ii\pi n(n-1)\a/2}\sum_{A\in\Ss_n} e^{2\ii\pi\a I(A)},$$
where $\Ss_n$ stands for the set of permutation matrices of size $n.$  Using the generating function of $I(A)$ for permutation matrices given e.g. in Corollary 3.5 of \cite{B}, further trigonometric simplifications entail 
$$\Dea^n(1,0+)\; =\; {\rm sf} (n-1)^2 V_n^\a,$$
as required by (\ref{Maindirect}).

\begin{Rq} {\em The formula (\ref{IKP}) also shows that
$$\Dea^n(x,x)\; =\; \frac{{\rm sf} (n-1)^2}{(2x\cos(\pi\a/2))^{n(n+1)}}\sum_{A\in\Aa_n} (4\sin^2(\pi\a/2))^{\mu(A)}\; > \; 0$$
for all $x > 0$ and $\a\in [0,1).$}
\end{Rq}

\section{Proof of part {\rm (c)}}

This last part of the theorem relies on Karlin's ETP criterion on the derivative determinant (see Theorem 2.2.6 in \cite{K}) which states that if $\Dea^k(x,y) > 0$ for every $k\in \{2, \ldots, n\}$ and $x,y >0,$ then $\Ka$ is ${\rm STP}_n.$ By an induction, we hence need  to show that
\begin{equation}
\label{Mainrev}
\a \le 1/n\wedge 1/(n^2-n-6)_+ \; \Rightarrow\; \Dea^n(x,y) > 0\quad\mbox{for all $x, y > 0$}.
\end{equation}
This will be obtained from (\ref{IKP}) and some considerations on ASM matrices, all to be found in \cite{B} and Section 2.1 of \cite{Be}. Again, for concision we will skip the dependence of the involved kernels in $(x,y)$. For any $A\in \Aa_n,$ introduce the statistics 
$$\nu(A)\; =\!\!\!\!\sum_{\tiny \begin{array}{c} 1\!\le\! i\!<\!i'\!\!\le\! n\\
1\!\le\! j\! <\!j'\!\!\le\! n\end{array}} \!\!\!\! \!\!\!\!A_{ij} A_{i'j'}\; =\; \frac{J(A)}{2}$$
(see pp. 5-6 in \cite{Be} for an explanation of the second equality). Set $\mu_n = \max\{\mu (A), \; A\in \Aa_n\}$ and notice that $\mu_n = (n-1)^2/4$ for $n$ odd and that $\mu_n = n(n-2)/4$ for $n$ even (see again 
 p. 6 in \cite{Be}). From (\ref{IKP}), it is clear that if 
\begin{equation}
\label{muAk}
F^\a_{n,k}\; =\; \sum_{\tiny \begin{array}{c}
A\!\in\!\Aa_n\\
\mu(A)\! =\! k\end{array}} \!\!\!\!   \Qa^{2\nu(A)}\bQa^{n(n-1) - 2\nu(A)-2k} \; > \; 0
\end{equation}
for all $k = 0\ldots \mu_n,$ then $\Dea^n$ is everywhere positive. Notice first that $F^\a_{n,k}$ is real since it writes 
$$\frac{1}{2}\!\!\sum_{\tiny \begin{array}{c}
A\!\in\!\Aa_n\\
\mu(A)\! =\! k\end{array}}\!\!\!\!  (\Qa^{2\nu(A)}\bQa^{n(n-1) - 2\nu(A)-2k} + \bQa^{2\nu(A)}\Qa^{n(n-1) - 2\nu(A)-2k}).$$
Indeed, setting $A^Q$ for the anticlockwise quarter-turn rotation of $A,$ one has $A^Q\in\Aa_n$ with $\mu(A^Q) = \mu(A)$ and $2\nu(A^Q) = n(n-1) - 2\nu(A) - 2\mu(A)$ (see p. 7 in \cite{Be}). \\

Suppose now that $k\neq 0.$ Then necessarily $n\ge 3, \,\nu(A)\ge 1$ and $\nu(A^Q)\ge 1$ (see p. 7 in \cite{Be}) and one gets the factorization
\begin{eqnarray*}
F^\a_{n,k} & = & \vert \Qa\vert^4 \!\!\!\!\sum_{\tiny \begin{array}{c}
A\!\in\!\Aa_n\\
\mu(A)\! =\! k\end{array}} \!\!\!\!\!\!  \Re\,(\Qa^{2(\nu(A)-1)}\bQa^{2(\nu(A^Q)-1)})\\
& = & \vert \Qa\vert^4
 \;\sum_{i=1}^{\nu_{n,k}} a^{i}_{n,k}\, \Re\,(\Qa^{2(i-1)}\bQa^{n(n-1)-2k - 2(i+1)})
\end{eqnarray*}
where $a^{i}_{n,k}$ are non-negative integer coefficients and $\nu_{n,k} = \max \{\nu(A), \; A\in \Aa_n, \mu(A) = k\}.$ Notice in passing that no closed expression for $a^{i}_{n,k}$ or even $\nu_{n,k}$ seem available in the literature. For $n =3,$ necessarily $k =1= \nu_{3,1}$ and there is only one summand, so that $F^\a_{3,1} = \vert \Qa\vert^4 > 0.$ For $n\ge 4$, it is sufficient to check that 
$$\Re\,(\Qa^{2(i-1)}\bQa^{n(n-1)-2k - 2(i+1)})\; > \; 0$$
for all $i=1\ldots \nu_{n,k}.$ Noticing that $n(n-1) -2k - 2(i+1)\le n^2-n-6$ 
and recalling that 
$$\Qa(x,y)\; =\; e^{\ii\pi\a/2} x + e^{-\ii\pi\a/2} y,$$ 
this inequality becomes clearly true when $\a\le 1/(n^2-n-6)$ after expanding the trigonometric polynomials, because all involved cosines are evaluated inside $[0,\pi/2].$  \\

It remains to consider the case $k =0.$ Reasoning exactly as above, for every $n\ge 2$ one obtains $F^\a_{n,0} > 0$ whenever $\a\le 1/n(n-1).$   The following proposition yields the optimal condition $\a\le 1/n$. Together with the above discussion, it concludes the proof of Part (c). \\

\begin{Prop} For every $n\ge 2,$ one has $F^\a_{n,0} > 0$ whenever $\a\le 1/n.$
\end{Prop}

\proof Using again the generating function of the inversion numbers of permutation matrices, we obtain the explicit formula
$$F^\a_{n,0}\; =\; \prod_{i=1}^n \lpa \frac{\bQa^{2i} - \Qa^{2i}}{\bQa^2 - \Qa^2}\rpa.$$
By an induction argument it is hence sufficient to prove that
$$\a\le 1/n\; \Rightarrow\; \frac{\bQa^{2n} - \Qa^{2n}}{\bQa^2 - \Qa^2}\; >\; 0$$
for every $n\ge 2.$ Setting $a = e^{2\ii \pi\a}$ and recalling the definition of $\Qa,$ we write
\begin{eqnarray*}
\frac{\bQa^{2n} - \Qa^{2n}}{\bQa^2 - \Qa^2} & = & \sum_{p+q =n-1}\sum_{j=0}^{2p} \sum_{k=0}^{2q} (a^{1/4}x)^j(a^{-1/4}y)^{2p-j} (a^{-1/4}x)^k(a^{1/4}y)^{2p-k}\\
& = & \sum_{i=0}^{2n-2} c_i\, x^i y^{2n-2-i}
\end{eqnarray*}
for some palindromic sequence $\{c_i, \, i= 0,\ldots, 2n-2\}$ which is given by
$$c_i\; =\; \sum_{\tiny\begin{array}{c}
p\!+\!q \!=\!n-1\\j\!+\!k \!=\!i\\j\!\le\! 2p\\k\!\le\! 2q\end{array}} a^{(j-k+q-p)/2}.$$
Summing appropriately and using some standard trigonometry, one can show that
$$c_i\; =\; \frac{\cos(\pi(n-1-i)\a)- \cos(\pi n\a)}{\sin^2(\pi\a)}$$
for every $i=0,\ldots, 2n-2.$ We omit the details. This completes the proof since the latter expressions are all positive whenever $\a\le 1/n.$ Notice also that these formul\ae \, show that the sequence $\{c_i\}$ is unimodal. 

\endproof

\section{Open questions}

In this section we give some heuristic reasons supporting the validity of the inclusion
\begin{equation}
\label{Mainreverse}
\a < 1/n \; \Rightarrow\; \Dea^n(x,y) > 0\quad\mbox{for all $x, y > 0$}
\end{equation}
which, by the ETP criterion and an immediate induction, is enough to show (\ref{Main}). First, we formulate a conjecture 
on the positivity of certain generating functions of ASM matrices which would entail (\ref{Mainreverse}), and we test it on the values $n=3,4,5.$ Second, we discuss more thoroughly Lascoux's factorization and settle two problems on the positivity on certain minors of rectangular matrices involving Chebyshev polynomials, whose solution would again entail (\ref{Mainreverse}).

\subsection{Alternating sign matrices} By (\ref{IKP}), the positivity of $\Dea^n$ amounts to that of 
$$\sum_{A\in\Aa_n} \Pa^{\mu(A)} \Qa^{2\nu(A)}\bQa^{2\nu(A^Q)},$$
with the notations of Section 4. The above sum can be expressed with the bivariate generating function
$$Z_n(x,y)\; =\; \sum_{A\in\Aa_n} x^{\nu(A)} y^{\mu(A)},$$
which is itself given as a certain functional determinant - see Formula (57) p. 21 in \cite{Be}. Unfortunately the latter determinant is in general very difficult to evaluate, except in certain particular cases. Its value at $x=y=1,$ counting the cardinal of $\Aa_n$, was the matter of a whole story which is told in the book \cite{B}. As in the proof of Part (c) let us now write
$$\sum_{A\in\Aa_n} \Pa^{\mu(A)} \Qa^{2\nu(A)}\bQa^{2\nu(A^Q)}\; =\; \sum_{k=0}^{\mu_n}\Pa^{\mu(A)}F^\a_{n,k},$$
recalling the notations $\mu_n = \max\{\mu (A), \; A\in \Aa_n\}$ and 
$$F^\a_{n,k}\; =\; \sum_{\tiny \begin{array}{c}
A\!\in\!\Aa_n\\
\mu(A)\! =\! k\end{array}} \!\!\!\!   \Qa^{2\nu(A)}\bQa^{2\nu(A^Q)}.$$
In view of Proposition 3, it is natural to raise the following question.

\begin{Conj} For every $n\ge 2$ and every $k = 0,\ldots,\mu_n,$ one has $F^\a_{n,k} > 0$ whenever $\a\le 1/n.$
\end{Conj}

The kernel $F^\a_{n,k}$ can be expressed in terms of the generating function
$$Z_{n,k}(x)\; =\; \frac{1}{k!} \frac{\partial^k Z_n}{\partial y^k}(x,0)$$
but unfortunately no tractable closed formula is known for the latter. Even its value at $x=1,$ which counts the number of elements of  $\Aa_n$ with $k$ negative entries, is known in closed form only in some cases - see Chapter 3 in \cite{LG} and the references therein. 

By Proposition 3, Conjecture 1 is true for $k = 0.$ In the preceding section, we also proved that $F^\a_{3,1} > 0$ whenever $\a\le 1/3.$ Let us now show the validity of Conjecture 1 for $n\ge 4$ and $k = \mu_n.$ Recall that $\mu_2 = 0$ and $\mu_3 = 1$ so that there is no loss of generality in considering $n\ge 4.$ 

\begin{Prop} For every $n\ge 4$ one has $F^\a_{n,\mu_n} > 0$ whenever $\a\le 1/n.$
\end{Prop}   

\proof First, suppose that $n = 2p+1$ is odd. Then $\mu_n = p^2$ and this concerns one single matrix with inversion number $\nu = p(p+1)/2$ (see  p. 6 in \cite{Be}). We deduce
$$Z_{n,\mu_n}(x)\; =\; x^{p(p+1)/2}\qquad\mbox{and}\qquad F^\a_{n,\mu_n}\; =\; \vert \Qa\vert^{2p(p+1)}\; >\; 0.$$
Second, suppose that $n = 2p$ is even. Then $\mu_n = p(p-1)$ and this concerns two matrices with inversion numbers $\nu = p(p+1)/2$ resp.  $\nu = p(p-1)/2$ (see again p. 6 in \cite{Be}). We deduce
$$Z_{n,\mu_n}(x)\; =\; x^{p(p-1)/2}(1+x^p)\qquad\mbox{and}\qquad F^\a_{n,\mu_n}\; =\; \vert \Qa\vert^{2p(p-1)}(\Qa^n +\bQa^n)\; >\; 0\quad \mbox{if $\a\le 1/n.$}$$
\endproof

Since there is no closed formula for $F^\a_{n,k}$ if $k\not\in\{0,\mu_n\},$ it seems quite difficult to establish Conjecture 1 for all $n,k.$ Let us conclude this paragraph in checking its validity for $n=4$ and $n=5.$ Since $\mu_4 =2$ and $\mu_5 =4$ we just have to consider the cases $ k=1$ resp. $k= 1,2,3.$

\medskip
$\bullet$ $Z_{4,1}(x)\, =\, 2x(1+x)^3$ and $F^\a_{4,1}\, =\, 2\vert \Qa\vert^4(\Qa^2 +\bQa^2)^3\, >\, 0$ if $\a\le 1/2.$

\medskip

$\bullet$ $Z_{5,1}(x) \, = \, x(3+14x + 35x^2+48x^3+48x^4 + 35 x^5 + 14x^6+ 3x^7),$ which rewrites 
$$3x\lpa\frac{x^5-1}{x^3-1}\rpa\,+\, 8x^2\lpa\frac{x^4-1}{x^2-1}\rpa\; +\; 10x^3\lpa\frac{x^4-1}{x^3-1}\rpa\; +\; 2x^4\lpa\frac{x^2-1}{x-1}\rpa$$
and yields
$$F^\a_{5,1}\, =\, 3\vert \Qa\vert^4\frac{F^\a_{5,0}}{F^\a_{3,0}}\;+\;8 \vert \Qa\vert^8 \frac{F^\a_{4,0}}{F^\a_{2,0}}  \;+\; 10\vert \Qa\vert^{12} \frac{F^\a_{4,0}}{F^\a_{3,0}}\;+\; 2\vert \Qa\vert^{16} F^\a_{2,0}.$$
The latter is positive for $\a\le 1/5,$ by Proposition 3.

\medskip

$\bullet$ $Z_{5,2}(x) \, = \, x(2+12x + 21x^2+24x^3+21x^4 + 12 x^5 + 2x^6),$ which rewrites 
$$
2x(1+x)(1+x+x^2)(1+x+x^2+x^3)\,+\, 6x^2(1+x^2)^2\, +\, 11x^3(1+x^2)$$
and yields 
$$F^\a_{5,2}\, =\, 2\vert \Qa\vert^4F^\a_{4,0}\,+\,6 \vert \Qa\vert^8 (\Qa^4 +\bQa^4)^2 \,+\, 11\vert \Qa\vert^{12} (\Qa^4 +\bQa^4).$$ 
The latter is positive for $\a\le 1/4,$ again by Proposition 3.

\medskip

$\bullet$ $Z_{5,3}(x)\, =\, x(1+6x^2+6x^3+x^5) \, =\, x(1+x)(3x^2 + (1+x^2)^2  - x(1-x)^2),$ so that
$$F^\a_{5,3}\; =\; \vert \Qa\vert^4(\Qa^2 +\bQa^2)(3 \vert \Qa\vert^8 + (\Qa^4 +\bQa^4)^2 - \vert \Qa\vert^4  (\Qa^2 -\bQa^2)^2 ).$$
Since $(\Qa^2 -\bQa^2)^2= -4(x-y)^2\sin^2(\pi\a),$ we see that $F^\a_{5,3} > 0$ if $\a\le 1/2.$ \\

Observe that the generating function $Z_{5,3}$ has vanishing terms so that $F^\a_{5,3}$ cannot be suitably written in terms of  $F^\a_{5,0},F^\a_{4,0},F^\a_{3,0},F^\a_{2,0},$ contrary to $F^\a_{5,2}$ and $F^\a_{5,1}.$ In general, testing the positivity of $F^\a_{n,k}$ seems to depend on bizarre rearrangements.

\subsection{Rectangular matrices} In this paragraph we study in more details the minors ${\rm A}_\sigma(1)$ and ${\rm B}_\sigma^\a(1),$ with the notations of Section 3.3. Indeed, formul\ae\, (\ref{separate}), (\ref{CaBi}) and (\ref{coeff}) show that the positivity of $\Dea^n$ is ensured by that of ${\rm A}_\sigma(1)$ and ${\rm B}_\sigma^\a(1).$ The analysis for ${\rm A}_\sigma(1)$ is easy.

\begin{Prop} For every $\sigma : \{1,\ldots,n\}\to\{1,\ldots, 2n-1\}$ increasing, one has ${\rm A}_\sigma(1) = {\rm A}_{\tilde \sigma}(1) > 0.$
\end{Prop}

\proof The generating function
$$\sum_{r\ge 0} C^r_{n-1+r} x^r \; =\; \frac{1}{(1-x)^n}$$
and Edrei's criterion - see Theorem 1.2 p. 394 in \cite{K} - show that the sequence $\{  C^r_{n-1+r}, \, r\ge 0\}$ is ${\rm TP}_\infty.$ In other words, the matrix $A_n(1)$ is ${\rm TP}_\infty$ viz. all its minors are non-negative. In particular, all coefficients ${\rm A}_\sigma(1)$ are non-negative. We now compute its exact positive value with the help of a formula of Gessel and Viennot on binomial determinants. Substracting the $i$-th row from the $(i-1)$-th row successively for $i = n\ldots 2$ and then repeating this operation for $i = (n-1)\ldots 2,$ $i= (n-2)\ldots 2,\,\ldots$ we see indeed from the Pascal relationships $C^{p+1}_r = C^{p+1}_{r+1} - C^p_r$ that the minor ${\rm A}_\sigma(1)$ is actually taken from the matrix
$$ \lcr C^{k-i}_{k-1}\rcr_{1\le i\le n, 1\le k\le 2n-1} \; =\;  \lcr C^{i}_{k}\rcr_{0\le i\le n-1, 0\le k\le 2n-2}.$$ 
The alternative formula mentioned at the bottom of p. 308 in \cite{GV} entails
$${\rm A}_\sigma(1)\; =\; \frac{\prod_{1\le i<j\le n} (\sigma_j -\sigma_i)}{{\rm sf} (n-1)}\cdot$$
In particular ${\rm A}_\sigma(1)= {\rm A}_{\tilde \sigma}(1)$ is always a positive integer (equal to one if and only if $\sigma_n -\sigma_1 = n-1).$ 

\endproof

The analysis for ${\rm B}_\sigma^\a(1)$ is however much more delicate. After transposition, we see that it is the $n\times n$ minor  obtained from the columns $\sigma_1, \ldots, \sigma_n$ in the horizontal matrix
$$ \lcr C^{k-i}_n U^\a_{2i-k}\rcr_{1\le i\le n, 1\le k\le 2n-1}.$$ 
Again, the generating function
$$\sum_{r\ge 0} C^r_n x^r \; =\; (1+x)^n$$
and Edrei's criterion show that the matrix $\lcr C^{k-i}_n \rcr_{1\le i\le n, 1\le k\le 2n-1}$ is ${\rm TP}_\infty.$ But since some $U^\a_{2i-k}$ are negative, nothing can be said a priori about the non-negativity of ${\rm B}_\sigma^\a(1).$ Transforming the rows $(L_1, \ldots, L_n)$ through the simultaneaous linear operations
$$L_i \,\rightarrow\, \sum_{j=i}^n C^j_n L_j$$
multiplies the $n\times n$ minors by a constant positive factor and yields a matrix whose $(i,k)$ coefficient is given by
$$\lacc\begin{array}{ll} \sum_{j=i}^n C^j_n C^{k-j}_n U^\a_{2j-k} & \mbox{if $1\le i \le n-1, i\le k\le n+i,$}\\
 C^{k-n}_n U^\a_{2n-k} & \mbox{if $ i = n, n\le k\le 2n-1,$}\\
0 & \mbox{otherwise.}
\end{array}\right.$$
Observe that the $(i,n+i)$ coefficient equals
$$ \sum_{j=i}^n C^j_n C^{n+i-j}_n U^\a_{2j-n-i}\; =\;  \sum_{j=i}^{[(n+i)/2]} C^j_n C^{n+i-j}_n (U^\a_{2j-n-i} + U^\a_{n+i-2j})\; =\; 0,$$
so that we have a band-matrix of width $n,$  whose $(i,k)$ coefficient is given by
$$\lacc\begin{array}{ll} \sum_{j=i}^n C^j_n C^{k-j}_n U^\a_{2j-k} & \mbox{if $1\le i \le n, i\le k\le n-1+i,$}\\
0 & \mbox{otherwise}.
\end{array}\right.$$
We set $B_{\a,n}$ for the above matrix, whose coefficients are all non-negative when $\a < 1/n.$ We also remark that $B_{\a,n}$ is persymmetric, viz. $B_{\a,n}(i,k) = B_{\a,n}(n+1-i, 2n-k)$ for any $(i,k).$ Indeed,
\begin{eqnarray*}
B_{\a,n}(i,k) - B_{\a,n}(n+1-i, 2n-k) & = &  \sum_{j=i}^n C^j_n C^{k-j}_n U^\a_{2j-k} \,+\,  \sum_{j=n+1-i}^n C^j_n C^{2n-k-j}_n U^\a_{2n-2j-k}\\
 & = &  \sum_{j=0}^n C^j_n C^{j+k-n}_n U^\a_{2n-2j-k}\; = \;  \sum_{j=n-k}^n C^j_n C^{j+k-n}_n U^\a_{2n-2j-k}
\end{eqnarray*}
and the last sum clearly vanishes. This persymmetry entails ${\rm B}_\sigma(1) = {\rm B}_{\tilde \sigma}(1).$ We now state a natural conjecture which would entail  (\ref{Mainreverse}).

\begin{Conj} {\em For every $n\ge 2,$ all coefficients  ${\rm B}_\sigma(1)$ are positive  whenever $\a<1/n.$}
\end{Conj}

The difficulty to prove this conjecture comes from the rectangular shape of the matrices $\{ B_{\a,n}\},$ which makes it seemingly impossible to use any kind of induction argument. Let us finally display the five first elements of the sequence $\{ B_{\a,n}\},$ whose common shape  reminds that of Trinidad and Tobago's national flag (the red background being the zeroes, the crucial white diagonal stripes Chebyshev polynomials of the second kind, and the black central parallelogram positive linear combinations thereof) and displays some spatial unimodality around the middle.
$$B_{\a,1}\; =\;(U^\a_1), \;B_{\a,2}\; =\; \lpa\begin{array}{ccc} 2U^\a_1 & U^\a_2 & 0 \\
0 &  U^\a_2 & 2 U^\a_1\end{array}\rpa,\; B_{\a,3}\; =\; \lpa\begin{array}{ccccc} 3U^\a_1 & 3 U^\a_2 & U^\a_3 & 0 & 0 \\
0 &  3U^\a_2 &  U^\a_3 + 9U^\a_1 & 3U^\a_2 & 0 \\
0 & 0 &  U^\a_3 & 3U^\a_2 & 3 U^\a_1\end{array}\rpa,$$
$$B_{\a,4}\; =\; \lpa\!\!\begin{array}{ccccccc} 4U^\a_1 & 6 U^\a_2 & 4U^\a_3 & U^\a_4 & 0 & 0& 0 \\
0 & 6U^\a_2  &  4U^\a_3 + 24U^\a_1 &  U^\a_4 + 16U^\a_2 & 4U^\a_3 & 0 & 0\\
0 & 0 & 4U^\a_3  &  U^\a_4 + 16U^\a_2 &  4U^\a_3 + 24U^\a_1 & 6U^\a_2 & 0 \\
0& 0 & 0 & U^\a_4 &  4U^\a_3 & 6U^\a_2 & 4U^\a_1\end{array}\!\!\rpa,$$
and 
$$B_{\a,5} =\lpa\!\!\begin{array}{ccccccc} 5U^\a_1 & 10 U^\a_2 & 10U^\a_3 & 5U^\a_4 & U^\a_5 & \cdots & 0 \\
0 & 10U^\a_2  &  10U^\a_3 + 50 U^\a_1&  5U^\a_4 + 50 U^\a_2 & U^\a_5 + 25U^\a_3 & \cdots & 0\\
0 & 0 & 10U^\a_3  &  5U^\a_4 + 50 U^\a_2 &  U^\a_5 + 25U^\a_3 + 100U^\a_1 & \cdots & 0\\
0 & 0 & 0& 5U^\a_4 &  U^\a_5 + 25U^\a_3 & \cdots & 0\\
0 & 0 & 0& 0 &  U^\a_5 & \cdots & 5U^\a_1\end{array}\!\!\rpa.$$
Using some elementary trigonometry, for these first five values of $n$ we could check that $B_{\a,n}$ is TP (viz. all its minors are non-negative) if $\a<1/n$  The following is hence natural

\begin{Conj} {\em For every $n\ge 1,$ the matrix $B_{\a,n}$ is  TP whenever $\a<1/n.$}
\end{Conj}

This conjecture is stronger than Conjecture 2 since it involves all minors. Notice that several criteria have appeared in recent years to prove the total positivity of a given matrix without checking every minor. See all the results mentioned in Section 2.5 of \cite{P} and also Theorem 2.16 therein for a criterion only on $2\times 2$ minors, interestingly also related to the zeroes of Chebyshev polynomials of the second kind. Unfortunately, none of this criteria seems particularly helpful in our situation.

\section{Proof of the corollary}

Setting $\fa(x) =\pa(1,x),$ we see from (\ref{PS}) that $\pa(t,x) = t^{-1/\a}\fa(xt^{-1/\a}),$ so that by Theorem 1.2.1. p.18 in \cite{K} the kernel $\pa(t,x)$ has the same sign-regularity over $(0,+\infty)\times (0,+\infty)$ than the kernel $\fa(e^{x-y})$ over $\rl\times\rl.$ In Paragraph 7.12.E p.390 of \cite{K} it is shown that 
$$\fa(e^{x-y})\,\in\, {\rm STP}_\infty\; \Leftrightarrow\; \fa(e^{x-y})\,\in\, {\rm TP}_\infty\; \Leftrightarrow\; \a\in\{1/2,1/3,\ldots, 1/n,\ldots\}$$
and the question is raised whether $\fa(e^{x-y})$ should be TP of some finite order when $\a$ is not the reciprocal of an integer. In \cite{ThS1}, we obtained the equivalence $\fa(e^{x-y})\,\in\, {\rm TP}_2\, \Leftrightarrow\, \a\le 1/2.$

Let us now prove the Corollary. We need to show that if $\a$ is not the reciprocal of an integer and $\a < 1/n,$ then $\fa(e^{x-y})\not\in {\rm SR}_n.$ If this were true, then the kernels $\fa(e^{y-x})$ and 
$$e^{x-y}\,\int_\rl  \fa(e^{x-u})\fa(e^{u-y})\, du$$
would also be ${\rm SR}_n$ by Theorem 1.2.1. and Lemma 3.1.1. in \cite{K}. However, a well-known fractional moment identification - see (3.1) in \cite{ThS1} and the references therein - shows that the latter kernel equals $\ga(x-y)$ with the notations of Section 2. Hence, we get a contradiction to Part (b) of the theorem.

\fin

We finish this paper with a natural conjecture on the total positivity of the positive stable kernel $\pa,$ which reformulates Karlin's question in a more precise manner. By Lemma 3.1.1. in \cite{K} and the same argument as in the proof of the corollary, this would also show the conjecture stated in the introduction. But we believe that this last conjecture is harder because the kernel $\pa$ is not explicit in general.

\begin{Conj} For every $n\ge 2$ one has 
$$\pa\,\in{\rm STP}_n\,\Leftrightarrow\,\pa\,\in{\rm SR}_n\,\Leftrightarrow\,\alpha\in \{1/2, 1/3, \ldots, 1/n,\ldots\}\;\; {\rm or}\;\; \alpha < 1/n.$$ 
\end{Conj}

\bigskip
\noindent
{\bf Acknowledgements.} I am grateful to Arno Kuijlaars for having lent me his personal copy of \cite{B}. Ce travail a b\'en\'efici\'e d'une aide de l'Agence Nationale de la Recherche portant la r\'ef\'erence ANR-09-BLAN-0084-01.

\end{document}